\documentclass[12pt]{amsart}
\usepackage{amsmath,amssymb}

\newtheorem{theorem}{Theorem}[section]

\newtheorem{corollary}[theorem]{Corollary}

\newtheorem{conjecture}[theorem]{Conjecture}

\theoremstyle{definition}

\theoremstyle{remark}

\numberwithin{equation}{section}
\numberwithin{theorem}{section}

\DeclareMathOperator{\vol}{vol}

\author{*Gang Tian}
\thanks{*The first author is supported in part by a National Science Foundation grant.}
\address{BICMR, Peking University, Yiheyuan Road 5, Beijing 100871, China}
\address{Department of Mathematics\\
Princeton University, Princeton NJ 08544, USA}
\email{gtian@math.pku.edu.cn}

\author{**Zhenlei Zhang}
\thanks{**The second author is supported by a grant of Beijing Municipal Commission of Education}
\address{School of Mathematics\\
Capital Normal University, Xisanhuan North Road 105, Beijing 100048, China}
\email{zhleigo@yahoo.com.cn}

\begin{document}

\title{Regularity of K\"ahler-Ricci flow}

\maketitle

\begin{abstract}
In this short note we announce a regularity theorem for K\"ahler-Ricci flow on a compact Fano manifold (K\"ahler manifold with positive first Chern class) and its application to the limiting behavior of K\"ahler-Ricci flow on Fano 3-manifolds. Moreover, we also present a partial $C^0$ estimate to the K\"ahler-Ricci flow under the regularity assumption, which extends previous works on K\"ahler-Einstein metrics and shrinking K\"ahler-Ricci solitons (cf. \cite{Ti90}, \cite{DoSu12}, \cite{Ti12}, \cite{PSS12}). The detailed proof will appear in \cite{TiZh13}.
\end{abstract}

\section{Introduction}

Let $M$ be a Fano $n$-manifold and $g_0$ be any K\"ahler metric with K\"ahler class $2\pi c_1(M)$. Consider the normalized K\"ahler-Ricci flow:
\begin{equation}\label{KRF}
\frac{\partial g}{\partial t}\,=\,g\,-\,{\rm Ric}(g),\hspace{0.5cm}g(0)\,=\,g_0.
\end{equation}
It was proved in \cite{Ca85} that \eqref{KRF} has a global solution $g(t)$ for $t\ge 0$. The main problem is to
understand the limit of $g(t)$ as $t$ tends to $\infty$.

By Perelman's non-collapsing result \cite{Pe02}, there is a uniform constant $\kappa=\kappa(g(0))>0$ satisfying:
\begin{equation}\label{volumenoncollap4}
\vol_{g(t)}(B_{g(t)}(x,r))\geq\kappa r^{2n},\hspace{0.5cm}\forall t\geq 0, r\leq 1.
\end{equation}
Since the volume stays the same along the K\"ahler-Ricci flow, this implies that for any sequence $t_i \rightarrow\infty$,
by taking a subsequence if necessary, $(M,g(t_i))$ converge to a limiting length space $(M_\infty,d)$ in the Gromov-Harsdorff topology:
\begin{equation}\label{e13}
(M,g(t_i))\stackrel{d_{GH}}{\longrightarrow}(M_\infty,d).
\end{equation}
The question is the regularity of $M_\infty$. A desirable picture is given in the following folklore conjecture \footnote{It is often referred as the Hamilton-Tian conjecture (see \cite{Ti97}).}.

\begin{conjecture}[\cite{Ti97}, also see \cite{ChTi06}]
\label{conj:HT}
$(M,g(t))$ converges (at least along a subsequence) to a shrinking K\"ahler-Ricci soliton with mild singularities.
\end{conjecture}

Here, "mild singularities" may be understood in two ways: (i) A singular set of codimension at least $4$, and (ii) a singular set of a normal variety. By extending the partial $C^0$-estimate conjecture \cite{Ti10} to
K\"ahler-Ricci flow, one can show that these two approaches are actually equivalent (see Section 3 or \cite{TiZh13}).

As pointed out by the second named author, this conjecture implies the Yau-Tian-Donaldson conjecture, in the case of Fano manifolds. The
conjecture states that if a Fano manifold $M$ admits a K\"ahler-Einstein metrics if it is K-stable. Recently, solutions were provided for this conjecture
in the case of Fano manifolds (\cite{Ti12}, also see \cite{CDS12a, CDS12b, CDS13}).


\section{Regularity of K\"ahler-Ricci flow}

Let $M$ be a Fano $n$-manifold and $g(t)$ a normalized
K\"ahler-Ricci flow in the K\"ahler class $2\pi c_1(M)$. Let
$(M_\infty,d)$ be a sequence limit of the K\"ahler-Ricci flow as
phrased in (\ref{e13}). The main regularity result is the following
theorem:

\begin{theorem}\label{regularity}
Suppose that for some uniform $p>n$ and $\Lambda< \infty$,
\begin{equation}\label{Lpricci4}
\int_M|Ric(g(t))|^pdv_{g(t)}\,\leq\,\Lambda.
\end{equation}
Then the limit $M_\infty$ is smooth outside a closed subset $\mathcal{S}$ of (real) codimension $\geq 4$ and $d$ is induced by a smooth K\"ahler-Ricci soliton $g_\infty$ on $M_\infty\backslash\mathcal{S}$. Moreover, $g(t_i)$ converge to $g_\infty$ in the $C^\infty$-topology outside $\mathcal{S}$.\footnote{The convergence with these properties is usually referred as the convergence in the Cheeger-Gromov topology, see \cite{Ti90} for instance.}
\end{theorem}

The proof of the theorem relies on Perelman's Pseudolocality theorem \cite{Pe02} of Ricci flow and a regularity theory for manifolds with $L^p$ bounded Ricci curvature ($p$ bigger than half dimension) and uniformly local volume noncollapsing condition as in (\ref{volumenoncollap4}). This is a generalization of the regularity theory of Cheeger-Colding \cite{ChCo96, ChCo97, ChCo00a} and Cheeger-Colding-Tian \cite{ChCoTi}. The proof can be carried out following the lines given in these papers under the framework established by Petersen-Wei \cite{PeWe97, PeWe00} on the geometry of manifolds with integral bounded Ricci curvature. Note that due to an example of Yang \cite{Ya92}, the uniformly volume noncollapsing condition (\ref{volumenoncollap4}) can not be replaced by a lower bound of total volume or local volume of metric balls of a definite size.

We shall show in \cite{TiZh13} that there is a uniform $L^4$ bound on the Ricci curvature along the K\"ahler-Ricci flow on any Fano manifold. Therefore, by the above regularity result, we have

\begin{corollary}
Conjecture \ref{conj:HT}, i.e., the Hamilton-Tian conjecture, holds for dimension $n\leq3$.
\end{corollary}

In the case of Del-Pezzo surfaces, Conjecture \ref{conj:HT} follows from \cite{tianzhu} and \cite{chenwang}.


\section{Partial $C^0$ estimate of K\"ahler-Ricci flow}

The partial $C^0$ estimate of K\"ahler-Einstein manifolds plays the key role in Tian's program to resolve the Yau-Tian-Donaldson conjecture, see \cite{Ti90}, \cite{Ti97}, \cite{Ti10}, \cite{DoSu12} and \cite{Ti12} for examples. An extension of the partial $C^0$ estimate to shrinking K\"ahler-Ricci solitons was given in \cite{PSS12}. These works are based on the compactness of Cheeger-Colding-Tian \cite{ChCoTi} and its generalizations to solitons by \cite{TiZh12}. We shall generalize these to the K\"ahler-Ricci flow on Fano manifolds in \cite{TiZh13} under the regularity assumption of the limit $M_\infty$.

Let $u(t)$ denote the Ricci potentials of the K\"ahler-Ricci flow $g(t)$ which satisfy
\begin{equation}
Ric(g(t))+\partial\bar{\partial}u(t)=g(t).
\end{equation}
The Hermitian metrics $\tilde{g}(t)=e^{-\frac{1}{n}u(t)}g(t)$ have $\omega(t)$, the K\"ahler forms of $g(t)$, as their Chern curvature forms. Let $H(t)$ be the induced metric on $K_M^{-l}$, the $l$-th power of the anti-canonical bundle ($l\geq 1$), and $D$ be the Chern connection of $H(t)$ on $K_M^{-l}$.

Let $\nabla$ and $\bar{\nabla}$ denote the $(1,0)$ and $(0,1)$ part of the Levi-Civita connection respectively. Then, at any time $t$, we have the Bochner type formula for $\sigma\in H^0(M,K_M^{-l})$
\begin{equation}\label{e18}
\triangle|\nabla\sigma|^2=|\nabla\nabla\sigma|^2+|\bar{\nabla}\nabla\sigma|^2
-\big((n+2)l-1\big)|\nabla\sigma|^2-\langle\partial\bar{\partial}u(\nabla\sigma,\cdot),\nabla\sigma\rangle
\end{equation}
and the Weitzenb\"{o}ch type formulas for $\xi\in C^\infty(M,T^{1,0}M\otimes K_M^{-l})$
\begin{equation}\label{e114}
\triangle_{\bar{\partial}}\xi=\bar{\nabla}^*\bar{\nabla}\xi+(l+1)\xi-\partial\bar{\partial}u(\xi,\cdot),
\end{equation}
\begin{equation}\label{e115}
\triangle_{\bar{\partial}}\xi=\nabla^*\nabla\xi-(n-1)l\xi,
\end{equation}
where $\triangle_{\bar{\partial}}$ is the Hodge Laplacian of $\bar{\partial}$. By \cite{ZhQ} there is a uniform bound of Sobolev constant under the K\"ahler-Ricci flow. Then apply the Moser iteration one can prove the gradient estimate to $\sigma\in H^0(M,K_M^{-l})$ and $L^2$ estimate to solutions $\bar{\partial}\vartheta=\xi\in C^\infty(M,T^{1,0}M\otimes K_M^{-l})$; see Lemmas 4.1 and 5.4 of \cite{Ti12}. The gradient estimate implies $\dim H^0(M,K_M^{-l})\leq N_l$ uniformly at any time $t$. We remark that Perelman's $C^1$ estimate to $u(t)$ \cite{SeTi08} will be used in the iteration arguments to cancel the bad terms containing $\partial\bar{\partial}u(t)$.

Now, let $\{s_{t,l,i}\}_{i=1}^{N_{t,l}}$ be an orthonormal basis of $H^0(M,K_M^{-l})$ with respect to the $L^2$ norm defined by $H(t)$ and Riemannian volume form, and put
\begin{equation}\label{e117}
\rho_{t,l}(x)=\sum_{i=1}^{N_{t,l}}|s_{t,l,i}|_{H}^2(x),\hspace{0.5cm}\forall x\in M.
\end{equation}
By using arguments similar to those in \cite{DoSu12} or \cite{Ti12}, we can prove

\begin{theorem}[Partial $C^0$ estimate]
If $(M,g(t_i))\stackrel{d_{GH}}{\longrightarrow}(M_\infty, g_\infty)$ as phrased in Theorem \ref{regularity}, then the partial $C^0$ estimate
\begin{equation}\label{e118}
\inf_{t_j}\inf_{x\in M}\rho_{t_j,l}(x)>0
\end{equation}
holds for a sequence of $l\rightarrow\infty$.
\end{theorem}

A direct corollary of this is to refine the regularity in Theorem \ref{regularity}.

\begin{theorem}\label{regularity-2}
Suppose $(M,g(t_i))\stackrel{d_{GH}}{\longrightarrow}(M_\infty, g_\infty)$ as phrased in Theorem \ref{regularity}. Then
$M_\infty$ is a normal projective variety and ${\mathcal S}$ is a subvariety of complex codimension at least $2$.
\end{theorem}

Finally, let us indicate how to deduce the Yau-Tian-Donaldson conjecture from the Hamilton-Tian conjecture.
Suppose $M$ is K-stable as defined in \cite{Ti97}. Then, under the K\"ahler-Ricci flow $g(t)$, we get a shrinking K\"ahler-Ricci soliton. This, together
with the uniqueness theorem on shrinking solitons, we can conclude that the Lie algebra of holomorphic vector fields on $M_\infty$ is reductive.
Then the K-stability implies the vanishing of Futaki invariant of $M_\infty$, consequently, the limit $(M_\infty,g_\infty)$ is K\"ahler-Einstein.
If $M_\infty$ is \textit{not} biholomorphic to $M$, then the eigenspaces of the first eigenvalues of $-\triangle_{g(t)}+g^{i\bar{j}}(t)\partial_iu(t)\partial_{\bar{j}}$ will converge to a subspace of potential functions on $M_\infty$ whose complex gradients are nontrivial holomorphic vector fields, cf. \cite{Zh11}. These vector fields induce the required degeneration of $M$ to $M_\infty$, with vanishing Futaki invariants. This gives a contradiction to the K-stability of $M$. So we have

\begin{theorem}
Suppose $M$ is K-stable. If $(M,g(t_i))\stackrel{d_{GH}}{\longrightarrow}(M_\infty, g_\infty)$ as phrased in Theorem \ref{regularity}, then $M$ coincides with $M_\infty$ and admits a K\"ahler-Einstein $g_\infty$.
\end{theorem}

In view of the regularity of low dimensional K\"ahler-Ricci flow in \S 2 we have

\begin{corollary}
The Yau-Tian-Donaldson conjecture holds for dimension $n\leq3$.
\end{corollary}

\end{document}